\documentstyle[12pt]{article}
\begin{document}

\newcommand{\mysection }[1]{\section{#1}\setcounter{equation}{0}}

\title{\bf The Existence of Type II Singularities for the Ricci Flow on $S^{n+1}$} \vspace{15mm}
\author{\bf}
\date{}
\maketitle \centerline{\large \bf Hui-Ling Gu and Xi-Ping Zhu}
\vspace{8mm} \centerline{{\large Department of Mathematics}}
\vspace{3mm} \centerline{{\large Sun Yat-Sen University }}
\vspace{3mm} \centerline{{\large Guangzhou, P. R. China}}
\vspace{3mm}
\vspace{15mm}

\noindent {\bf Abstract } In this paper we prove the existence of
Type II singularities for the Ricci flow on $S^{n+1}$ for all $n
\geq 2$.

\begin{center}

{\bf \large \bf{1. Introduction}}
\end{center}

 In this paper, we consider the Ricci flow
$$
      \left\{
       \begin{array}{lll}
\frac{\partial g_{ij}}{\partial t}=-2R_{ij},
          \\[4mm]
  g_{ij}(0)=\hat g_{ij},
       \end{array}
    \right.
  \eqno(1.1)
$$
starting from a given compact Riemannian manifold (${M},\hat g$).
This is a nonlinear (degenerate) parabolic system on metrics. In
the seminal paper \cite{Ha82}, Hamilton proved the Ricci flow
admits a unique solution on a maximal time interval $[0,T)$ so
that either $T=+\infty$ or $T<+\infty$ and $|Rm|$ is unbounded as
$t\rightarrow T$. We call such a solution $g(t)$ a maximal
solution of the  Ricci flow. If $T<+\infty$ and the curvature
becomes unbounded as $t$ tends to $T$, we say the maximal solution
develops singularities as $t$ tends to $T$ and $T$ is a singular
time. It is well-known that the Ricci flow generally develops
singularity.

If a solution $(M,g(t))$ to the Ricci flow develops singularities
at a maximal time $T < +\infty$, according to Hamilton
\cite{Ha95F}, we say it develops a \mbox{\textbf{Type I
singularity}} if
$$ \ \sup_{t\in [0,T)}(T-t)K_{max}(t)<+\infty,\
$$
and say it develops a  \mbox{\textbf{Type II singularity}} if
$$ \ \sup_{t\in [0,T)}(T-t)K_{max}(t)=+\infty,
$$
where $K_{max}(t) = \max \{|Rm(x,t)| \ | \ x \in M \}$.

Clearly, a round sphere, or more generally a finite product of
several space-forms with positive curvature, shrinks to form Type I
singularities. In \cite{Ha88, Chow}, Hamilton and Chow proved the
Ricci flow on two-sphere $S^2$ (with an arbitrary metric) always
develops a Type I singularity and shrinks to a round point. In
\cite{Ha82, Ha86}, Hamilton proved the Ricci flow on a compact
three-manifold with positive Ricci curvature, or a compact
four-manifold with positive curvature operator, develops a Type I
singularity and shrinks to a round point; recently,
B$\ddot{o}$hm-Wilking \cite{BW} had shown that the Ricci flow on a
general compact $n$-dimensional Riemannian manifold with positive
curvature operator also develops a Type I singularity and shrinks to
a round point.

Intuitively, a compact manifold with the shape like a dumbbell will
develop a Type I singularity in the neck part. In views of the work
\cite{Ha97} of Hamilton on four-manifolds with positive isotropic
curvature (see also \cite{CZ05}), a Type I singularity with
neckpinch is expected. Indeed, such examples of Type I singularities
with neckpinch for the mean curvature flow were known more than
fifteen years ago (see for example \cite{Gr} and \cite{An-m}). It is
very surprising that the existence of Type I singularities with
neckpinch for the Ricci flow was only known very recently. The first
rigorous examples of Type I singularity with neckpinch for the Ricci
flow were constructed by Miles Simon \cite{MS} on \emph{noncompact}
warped products $R\times_fS^n$. In \cite{FIK}, Feldman-Ilmanen-Knopf
constructed another family of rigorous examples of Type I with
neckpinch on the holomorpic line bundle $L^{-k}$ over $CP^{n-1}$.
Both of these families of examples live on \emph{noncompact}
manifolds. For the Type I singularity with neckpinch on compact
manifolds, the first rigorous examples were given by Sigurd Angenent
and Dan Knopf in \cite{AK} by constructing suitable rotationally
symmetric metrics on $S^{n+1}$, where the definition of a
rotationally symmetric metric is the following:
 \vskip
0.3cm \noindent {\bf Definition 1.1}   A metric $g$ on $I\times
S^n$, where $I$ is an interval, is called \textbf{rotationally
symmetric} if it has the following form:$$g=\varphi (x)^2dx^2+\psi
(x)^2g_{can},\qquad x\in I,$$ where $g_{can}$ is the standard
metric of the round sphere $S^n$ with constant (sectional)
curvature 1.

For the Type II singularity for the Ricci flow, a rigorous example
on $R^2$ was recently given by Daskalopoulos and Hamilton in
\cite{DH}. However, no rigorous examples of Type II singularity for
the Ricci flow on compact manifolds have yet appeared. We remark
that some beautiful intuitions of the forming of Type II singularity
were described and explained by Chow-Knopf in \cite{CK} and Topping
in \cite{T}. (For the mean curvature flow, the existence of Type II
singularities was already justified by Altschuler-Angenent-Giga
\cite{AAG} and Angenent-Vel$\acute{a}$zquez \cite{An-c}.)

The purpose of this paper is to demonstrate the existence of Type II
singularity on compact manifolds, in particular for
rotationally-symmetric initial metrics on $S^{n+1}$. Our main result
is the following:

\vskip 0.3cm \noindent {\bf Theorem 1.2} \emph{ For each $n\geq 2$,
there exist rotationally-symmetric metrics on $S^{n+1}$ such that
the Ricci flow starting at the metrics develop Type II singularities
at some times $T<+\infty$.}\vskip 0.3cm

This paper contains four sections. In Section 2, we recall some
useful estimates of Angenent-Knopf \cite{AK} on rotationally
symmetric solutions to the Ricci flow. In general, to understand
the structure of singularities, one usually needs to get a
classification for gradient shrinking solitons. The recent work
\cite{P2} of Perelman gives a complete classification to
positively curved gradient shrinking Ricci soliton in dimension
three. In Section 3 we will extend Perelman's classification to
higher dimensions in the class of rotationally symmetric metrics.
Finally in Section 4, based on the generalized classification, we
will prove the main result Theorem 1.2.

Our work in this paper benefits from a conversation with Professor
R. S. Hamilton, who suggested the second author to consider the
class of rotationally symmetric metrics. The second author is
partially supported by NSFC 10428102 and NKBRPC 2006CB805905.

\begin{center}{ \bf \large \bf {2. Angenent-Knopf's Estimates}}\end{center}

Consider a rotationally-symmetric metric
$$g=\varphi (x)^2dx^2+\psi (x)^2g_{can} \eqno(2.1)$$
on the set $(-1,1)\times{S^n}$, in which $g_{can}$ is the metric
of constant sectional curvature 1 on $S^n$. The coordinate $x$ is
ungeometric, a more geometric quantity is the distance $s$ to the
equator given by
$$s(x)=\int_0^x\varphi(x)dx.$$ Then
$$\frac{\partial}{\partial s}=\frac{1}{\varphi(x)}\frac{\partial}{\partial
x}$$ and $$ds=\varphi(x)dx.$$ With this notation the metric is
$$g=ds^2+\psi^2g_{can}. \eqno(2.2)$$
 In order to extend $g$ to be
a smooth Riemannian metric on $S^{n+1}$, it is sufficient  and
necessary to impose the boundary conditions: $$\psi(\pm1)=0,\qquad
\lim_{x\rightarrow\pm1}\psi_{s}(x)=\mp1$$ and
$$\lim_{x\rightarrow\pm1}\frac{d^{2k}\psi(x)}{ds^{2k}}=0$$ for all $k = 1, 2, \cdots$. The
Riemannian curvature tensor of (2.2) is determined by the
sectional curvatures
$$K_0=-\frac{\psi_{ss}}{\psi}$$ of the $n$ 2-planes perpendicular
to the spheres $\{x\}\times S^n$, and the sectional curvatures
$$K_1=\frac{1-\psi_s^2}{\psi^2}$$ of the $\frac{n(n-1)}{2}$
2-planes tangential to these spheres. In the ungeometric
coordinate $x$ the Ricci tensor of the metric $g$ given by (2.1)
is
$$Ric=n\{-\frac{\psi_{xx}}{\psi}+\frac{\varphi_x\psi_x}{\varphi\psi}\}(dx)^2+
\{-\frac{\psi\psi_{xx}}{\varphi^2}-\frac{(n-1)\psi_x^2}{\varphi^2}+\frac{\psi\varphi_x\psi_x}{\varphi^3}
+n-1\}g_{can}.$$ In the geometric coordinate this simplifies to
$$Ric=(nK_0)ds^2+\psi^2[K_0+(n-1)K_1]g_{can}.$$The scalar
curvature is given by $$R=2nK_0+n(n-1)K_1.$$ The above
computations can be found in \cite{AK} or the textbook \cite{PP}.

Suppose we have a time dependent family of metrics $g(\cdot,t)$
having the form (2.1). Then the family $g(\cdot,t)$ satisfies the
Ricci flow if and only if $\varphi$ and $\psi$ evolve by:
$$\varphi_t=n\frac{\psi_{ss}}{\psi}\varphi, \eqno(2.3)$$
$$\psi_t=\psi_{ss}-(n-1)\frac{1-\psi_s^2}{\psi}. \eqno(2.4)$$








According to Angenent-Knopf \cite{AK}, the (interior) local
minimal points of the function $x\mapsto \psi(x,t)$ are called
``\textbf{necks}" and the (interior) local maximal points are
called ``\textbf{bumps}". As long as the solution exists at a time
$t$, the radius of the smallest neck is given by
$$r_{min}(t)=\min\{\psi(x,t)|\psi_x(x,t)=0\}.$$
Of course, if the solution has no necks at the time $t$, we let
$r_{min}$ not be defined. Denote by $x_+(t),x_-(t)$ the right-most
bump (i.e. the largest local maximal point on $(-1,+1)$) and
left-most bump (i.e. the least local maximal point on $(-1,+1)$)
respectively. The region right of $x_+(t)$ and left of $x_-(t)$ are
called the ``\textbf{right polar cap}" and ``\textbf{left polar
cap}" respectively. In \cite{AK},  Angenent and  Knopf obtained
several useful estimates for the Ricci flow via the equations (2.3)
and (2.4). We recall some of their estimates as follows.

\vskip 0.2cm\noindent {\bf Proposition 2.1} (Angenent-Knopf
\cite{AK}) \emph{ Let $g(t)$ be a solution to the Ricci flow of the
form (2.2) such that $|\psi_s|\leq 1$ and the scalar curvature $R>0$
and $\psi_s$ has finitely many zeroes initially.  Then}


(1) (Proposition 5.1 of \cite{AK}) \emph{ As long as the solution
exists, $|\psi_s|\leq1$.}

(2) (Lemma 7.1 of \cite{AK}) \emph{ There exists $C=C(n,g(0))$ such
that as long as the solution exists,$$|Rm|\leq\frac{C}{\psi^2}.$$}


(3) (Lemma 5.6 and Lemma 7.2 of \cite{AK}) \emph{ If the left polar
cap is strictly concave (i.e., $\psi_{ss}<0$)at initial, then as
long as the solution exists, left polar cap exists and remains
strictly concave, and $D=\lim_{t\nearrow T}\psi(x_-(t),t)$ exists.
Furthermore, no singularity occurs on the left polar cap if $D >
0$.}


(4) (Lemma 9.1 of \cite{AK}) \emph{ There exists $C=C(n,g_0)$ such
that
$$\frac{K}{L}[logL+2-logL_{min}(0)]\leq C,$$ where $K=-K_0=\frac{\psi_{ss}}{\psi}$ and
$L=K_1=\frac{1-\psi_s^2}{\psi^2}$.}$$\eqno \#$$


\begin{center}{ \bf \large \bf{ 3. Classification of Shrinking Solitons}}\end{center}

To understand the structure of singularities, one usually needs to
get a classification for gradient shrinking solitons. In \cite{P2},
Perelman obtained a complete classification for nonnegatively curved
gradient shrinking soliton in dimension three. An open question is
how to generalize Perelman's classification to higher dimensions. In
the next proposition, we obtain such a classification for the class
of rotationally symmetric solitons.

\vskip 0.2cm\noindent {\bf Proposition 3.1} \emph{ Let $(M,
g_{ij}(t))$, $-\infty <t<0$, be a nonflat gradient shrinking
soliton to the Ricci flow on a complete $(n+1)$-dimensional
manifold and assume the metric $g_{ij}(t)$ is rotationally
symmetric. Suppose $(M,g_{ij}(t))$ has bounded and nonnegative
sectional curvature and is $\kappa$-noncollapsed on all scales for
some $\kappa>0$. Then $(M,g_{ij}(t))$ is one of the followings: }

(i) \emph{the round sphere $S^{n+1}$;}

(ii) \emph{the round infinite cylinder $(-\infty,+\infty) \times
S^n $.}

\vskip 0.1cm \noindent{\bf Proof.} Note that for a rotationally
symmetric metric, the nonnegativity of sectional curvatures is
equivalent to the nonnegativity of curvature operator. Indeed, we
can choose a coordinate system $(x^0,x^1,\cdot\cdot\cdot,x^n)$
(where $x^0$ is the radial direction and $x^i, i = 1,
\cdot\cdot\cdot, n,$ are the spherical directions) on $M$ such that
all components of the Riemannian curvature tensor vanish in the
coordinate system except the sectional curvatures
$R_{i0i0}=\psi^2K_0$ and $R_{ijij}=\psi^4K_1$ $ (i\neq j)$, and then
the equivalence follows directly from Proposition 1.1 and 1.2 of
\cite{PP}.

 Firstly, we
consider the case that the gradient shrinking soliton is compact
and has strictly positive sectional curvature everywhere. By the
Theorem 1 in \cite{BW} we see that the compact gradient shrinking
soliton is getting round and tends to a space form (with positive
constant curvature) as the time tends the maximal time $t=0$.
Since the shape of a gradient shrinking Ricci soliton is
unchanging up to reparameterizations and homothetical scalings,
the gradient shrinking soliton has to be the round $(n+1)$-sphere
$S^{n+1}$.

Next, we consider the case that the sectional curvature of the
nonflat gradient shrinking soliton vanishes somewhere. Note that a
rotationally symmetric metric is defined on $I \times S^n$ for some
interval $I$. By Hamilton's strong maximum principle in \cite{Ha86},
we know that the soliton splits off a line and then the soliton is
the round cylinder $R \times S^n $. (We remark that $R^k \times
S^{n+1-k}$ is not rotationally symmetric if $k>1$.)

Finally we want to exclude the case that the gradient shrinking
soliton is noncompact and has strictly positive sectional
curvature everywhere.

Suppose there is a complete $(n+1)$-dimensional noncompact
$\kappa$-noncollapsed gradient shrinking soliton $g_{ij}(t)$,
$-\infty<t<0$, satisfies $$\nabla_i\nabla_j
f+R_{ij}+\frac{1}{2t}g_{ij}=0,\ \ \mbox{on}\ -\infty<t<0, \eqno
(3.1)$$ everywhere for some function $f$ and
$g(t)=ds^2+\psi^2(s,t)g_{can}$ and with bounded and positive
sectional curvature at each time $t\in(-\infty,0)$. Let us
consider the shrinking soliton at the time $t=-1$. Arbitrarily fix
a point $x_0$ in $M$. By the same arguments as in the proof of
Lemma 1.2 of Perelman \cite{P2} (or see the proof of Lemma 6.4.1
of \cite{CZ} for the details), one has the followings:

(1) at large distance from the fixed point $x_0$ the function $f$
has no critical point, and its gradient makes small angle with the
gradient of the distance function from $x_0$;

(2) at large distance from $x_0$, the scalar curvature $R$ is
strictly increasing along the gradient curves of $f$, and
$$\limsup_{d_{(-1)}(x,x_0)\rightarrow+\infty}R(x,-1) \leq \frac{n}{2};$$

(3) the volume of the level set of $f$
 satisfies $$Vol\{f=a\}<Vol(S^n(\sqrt{2(n-1)}))\eqno(3.2)$$ for all large
enough $a$.

In the three-dimension case, Perelman (in Lemma 1.2 of \cite{P2})
argued by using Gauss-Bonnet formula to the level set $\{f=a\}$ to
derive a contradiction. But now we are considering the general
dimensional case, in particular, the (generalized) Gauss-Bonnet
formulas are not available. So we need a new argument in the
following.

By using Gauss equation and (3.1), the intrinsic sectional
curvature $\tilde{R}_{ijij}$ of the level set $\{f=a\}$ can be
compute as
$$\arraycolsep=1.5pt\begin{array}{rcl}
\tilde{R}_{ijij}&=&R_{ijij}+(h_{ii}h_{jj}-h_{ij}^2)\\[4mm]
&=&R_{ijij}+\frac{1}{|\nabla f|^2}(f_{ii}f_{jj}-f_{ij}^2)\\[4mm]
&\leq& R_{ijij}+\frac{1}{4|\nabla f|^2}(f_{ii}+f_{jj})^2\\[4mm]
&=&R_{ijij}+\frac{1}{4|\nabla f|^2}(1-R_{ii}-R_{jj})^2.
\end{array}\eqno(3.3)$$

Denote by $X=\frac{\nabla f}{|\nabla f|}$ the unit normal vector
to the level set $\{f=a\}$. Then set
$$X=\delta^0\frac{\partial}{\partial x^0}+\delta^\alpha\frac{\partial}{\partial
x^\alpha},$$ and $$e_i=u_i^0\frac{\partial}{\partial
x^0}+u_i^\alpha\frac{\partial}{\partial x^\alpha},\qquad
i=1,2,\cdot\cdot\cdot,n.$$ where the summation convention of summing
over repeated indices is used and $\{x^0,x^1,\cdot\cdot\cdot,x^n\}$
is the local coordinate on the $(n+1)$-dimensional rotationally
symmetric manifold with $g=ds^2+\psi^2g_{can}$ with $x^0=s\in R$ and
$P=(x^1,\cdot\cdot\cdot,x^n)\in S^n$ and
$g_{\alpha\beta}=\delta_{\alpha\beta}$ at $(s,P)$. In these
coordinates all components of the Riemann tensor and Ricci tensor
vanish except $R_{\alpha 0\alpha 0}=K_0$ and
$R_{\alpha\beta\alpha\beta}=K_1 (\alpha\neq\beta)$  and
$R_{00}=nK_{0}$ and $R_{\alpha\alpha}=K_0+(n-1)K_1,$ ($\alpha
=1,2,\cdot\cdot\cdot,n$) where $K_0=-\frac{\psi_{ss}}{\psi}$ and
$K_1=\frac{1-\psi_s^2}{\psi^2}$. And then the scalar curvature
$R=2nK_0+n(n-1)K_1$. So we have
$$\arraycolsep=1.5pt\begin{array}{rcl}
R_{ijij}&=&R(u_i^\alpha \frac{\partial}{\partial
x^\alpha},u_j^\beta\frac{\partial}{\partial
x^\beta},u_i^\gamma\frac{\partial}{\partial x^\gamma},
u_j^\eta\frac{\partial}{\partial x^\eta})\\[4mm]
&=&\sum_{\alpha\beta\gamma\eta}u_i^\alpha u_j^\beta u_i^\gamma
u_j^\eta R_{\alpha\beta\gamma\eta}\\[4mm]
&=&\sum_{\alpha\beta}(u_i^\alpha u_j^\beta )^2 R_
{\alpha\beta\alpha\beta}-\sum_{\alpha\beta}u_i^\alpha u_j^\beta
u_i^\beta u_j^\alpha R_{\alpha\beta\alpha\beta}\\[4mm]
&=&\sum_{\beta=1}^n[(u_i^0)^2(u_j^\beta)^2+(u_i^\beta)^2(u_j^0)^2]K_0+\sum_{\alpha,\beta=1}
^n(u_i^\alpha u_j^\beta)^2K_1\\[4mm]
&&-2\sum_{\beta=1}^n u_i^0u_j^0u_i^\beta u_j^\beta
K_0-\sum_{\alpha,\beta=1}^n u_i^\alpha u_j^\beta u_i^\beta
u_j^\alpha K_1\\[4mm]
&=&(\sum_{\beta=1}^n[(u_i^0)^2(u_j^\beta)^2+(u_i^\beta)^2(u_j^0)^2]+2(u_i^0)^2(u_j^0)^2)K_0\\[4mm]
&&+\sum_{\alpha,\beta=1}^n [(u_i^\alpha u_j^\beta)^2-u_i^\alpha
u_j^\beta u_i^\beta u_j^\alpha]K_1\\[4mm]
&=&[(u_i^0)^2+(u_j^0)^2]K_0+\sum_{\alpha,\beta=1}^n[(u_i^\alpha
u_j^\beta)^2-u_i^\alpha u_j^\beta u_i^\beta u_j^\alpha]K_1\\[4mm]
&=&[(u_i^0)^2+(u_j^0)^2]K_0+[(1-(u_i^0)^2)(1-(u_j^0)^2)-(u_i^0)^2(u_j^0)^2]K_1\\[4mm]
&=&[(u_i^0)^2+(u_j^0)^2]K_0+[1-(u_i^0)^2-(u_j^0)^2]K_1.
\end{array}\eqno(3.4)$$ where in the fifth and sixth equalities we
used $$\sum_{\beta=1}^n u_i^\beta u_j^\beta=-u_i^0 u_j^0$$ and
$$\sum_{\beta=0}^n (u_j^\beta)^2=1$$
since $\{e_1,\cdots,e_n\}$ is an orthonormal basis of the level set
$\{f=a\}$. Then by (3.4) we have
$$\arraycolsep=1.5pt\begin{array}{rcl}
R_{ii}&=&R_{iXiX}+\sum_{j=1,j\neq i}^n R_{ijij}\\[4mm]
&=&[(u_i^0)^2+(\delta^0)^2]K_0+[1-(u_i^0)^2-(\delta^0)^2]K_1\\[4mm]
&&+\sum_{j=1,j\neq i}^n [(u_i^0)^2+(u_j^0)^2]K_0+\sum_{j=1,j\neq
i}^n [1-(u_i^0)^2-(u_j^0)^2]K_1\\[4mm]
&=&[n(u_i^0)^2+1-(u_i^0)^2]K_0+[n-n(u_i^0)^2-1+(u_i^0)^2]K_1\\[4mm]
&=&[1+(n-1)\varepsilon]K_0+[(n-1)(1-\varepsilon)K_1]
\end{array}\eqno(3.5)$$ where $\varepsilon=(u_i^0)^2\ll 1$, if $a$
is large enough.

Obviously by (3.5) we get
$$R_{ii}<2K_0+(n-1)K_1=\frac{R}{n}<\frac{1}{2},$$ and then
$$1-R_{ii}-R_{jj}>0.$$ Again by (3.5) we know $$R_{ii}>K_0+(n-1)(1-\varepsilon)
K_1\eqno(3.6)$$ and by (3.4) we know $$R_{ijij}<2\varepsilon
K_0+(1-\varepsilon)K_1.\eqno(3.7)$$ Hence by (3.3), (3.6) and (3.7)
$$\arraycolsep=1.5pt\begin{array}{rcl}
\tilde{R}_{ijij}&\leq& R_{ijij}+\frac{1}{4|\nabla
f|^2}(1-R_{ii}-R_{jj})^2\\[4mm]
&<&2\varepsilon K_0+(1-\varepsilon)K_1+\frac{1}{4|\nabla
f|^2}[1-2(K_0+(n-1)(1-\varepsilon)K_1)]^2\\[4mm]
&=&2\varepsilon
K_0+\frac{1}{2(n-1)}\{1-2(K_0+(n-1)(1-\varepsilon)K_1+2(n-1)(1-\varepsilon)K_1\\[4mm]
&&-[1-2(K_0+(n-1)(1-\varepsilon)K_1)]\}+\frac{1}{4|\nabla
f|^2}[1-2(K_0+(n-1)(1-\varepsilon)K_1)]^2\\[4mm]
\end{array}$$
$$\arraycolsep=1.5pt\begin{array}{rcl}
&=&\frac{1}{2(n-1)}\{1-2(1-2\varepsilon
(n-1))K_0-[1-2(K_0+(n-1)(1-\varepsilon)K_1)]\\[4mm]
&&+\frac{n-1}{2|\nabla f|^2}[1-2(K_0+(n-1)(1-\varepsilon)K_1)]^2\}\\[4mm]
&<&\frac{1}{2(n-1)}
\end{array}\eqno(3.8)$$ for sufficiently large $a$, since $2(1-2\varepsilon
(n-1))K_0>0$ and $1-2(K_0+(n-1)(1-\varepsilon)K_1)>0$ and $|\nabla
f|$ is large as $a$ large. Then by (3.8) and the volume comparison
theorem we know $$Vol\{f=a\}>Vol(S^n(\sqrt{2(n-1)}))$$ for large
enough $a$ and then it is a contradiction with (3.2).

Therefore we have proved the proposition.$$\eqno \#$$

\begin{center}{ \bf \large \bf{ 4. Type II Singularity Happens}}\end{center}

Suppose we have a family of rotationally symmetric solutions
$$\{(S^{n+1},g_\alpha(t))|  \alpha\in[0,1]\}$$ of the Ricci flow
with $g_\alpha(0)=ds^2+\psi_\alpha^2g_{can}$, $\alpha \in [0,1]$,
where $g_{can}$ is the standard metric of constant sectional
curvature 1 on $S^{n}$. We specify the initial metrics as follows.

When $\alpha=1$, let the initial metric $g_1(0)$ be a symmetric
dumbbell with two equally-sized hemispherical regions joined by a
thin neck. By the work in \cite{AK}, we can assume the two
hemispheres are suitably large and the neck is suitably thin so that
this initial metric $g_1(0)$ leads to a neckpinch singularity of the
Ricci flow at some time $T_1<+\infty$. (see Figure 1.)
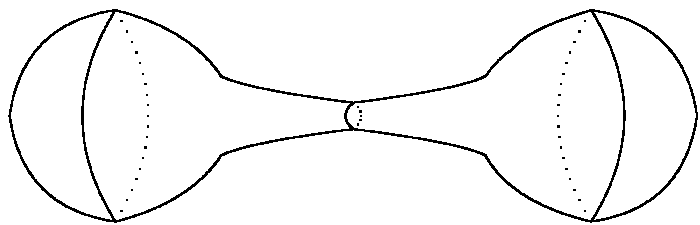
\begin{figure}[h]
\begin{picture}(350,90)
\qbezier(50,50)(55,85)(90,90)
\qbezier(50,50)(55,15)(90,10)\qbezier(90,90)(65,50)(90,10)\qbezier[20](90,90)(115,50)(90,10)
\qbezier(90,90)(118,83)(130,65)
\qbezier(90,10)(118,17)(130,35)\qbezier(180,55)(174,50)(180,45)\qbezier[6](180,55)(186,50)(180,45)
\qbezier(130,65)(140,60)(180,55) \qbezier(130,35)(140,40)(180,45)
\qbezier(180,55)(220,60)(230,65) \qbezier(180,45)(220,40)(230,35)
\qbezier(240,75)(246,83)(270,90)\qbezier(230,65)(235,72)(240,75)
\qbezier(230,35)(240,17)(270,10)\qbezier(270,90)(295,50)(270,10)\qbezier[20](270,90)(245,50)(270,10)
\qbezier(270,90)(305,85)(310,50) \qbezier(270,10)(305,15)(310,50)
\end{picture}
\caption{A neckpinch forming}
\end{figure}

When $\alpha=0$, let the initial metric $g_0(0)$ be a lopsided and
degenerate dumbbell where $g_0(0)=ds^2+\psi_0^2g_{can}$ with
$\psi_0(0)$ has only one bump and it is nonincreasing on the right
polar cap and strictly concave on the left polar cap. (see Figure
2.)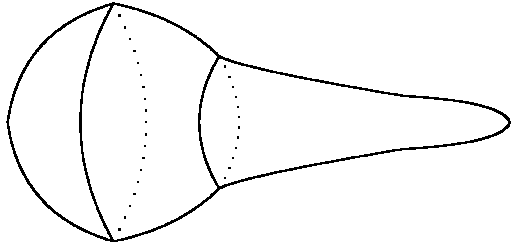
\begin{figure}[h]
\begin{picture}(350,105)
\qbezier(50,60)(55,95)(90,105)
\qbezier(50,60)(55,25)(90,15)\qbezier(90,105)(65,60)(90,15)\qbezier[20](90,105)(115,60)(90,15)
\qbezier(90,105)(115,100)(130,85)\qbezier(90,15)(115,20)(130,35)\qbezier(130,85)(115,60)(130,35)
\qbezier(130,85)(140,80)(200,70)
\qbezier(130,35)(140,40)(200,50)\qbezier[13](130,85)(145,60)(130,35)
\qbezier(200,70)(237,68)(240,60) \qbezier(200,50)(237,52)(240,60)
\end{picture} \caption{A degenerate dumbell}
\end{figure}

Clearly, we may choose the $g_1(0), g_0(0)$ to have positive scalar
curvatures. Let $\{ g_\alpha(0) \ | \ \alpha \in [0,1] \}$ (see
Figure 3.) be a smooth family of dumbbells (including degenerate
dumbbells) connecting the $g_1(0)$ to the $g_0(0)$ and satisfying
the followings: \vskip 0.1cm\noindent (i)  for each $\alpha \in
[0,1]$, $\psi_\alpha(0)$ has exactly two bumps or one bump, \vskip
0.1cm\noindent (ii)  for each $\alpha \in [0,1]$,
$(\psi_\alpha)_s(0)$ has only finitely many zeros, and satisfies
$$|(\psi_\alpha)_s|(0) \leq 1,$$ \vskip 0.1cm\noindent(iii) for each
$\alpha \in [0,1]$, $\psi_\alpha(0)$ is strictly concave on the left
polar cap, \vskip 0.1cm\noindent(iv) each initial metric
$g_\alpha(0)$, $\alpha \in [0,1]$, has positive scalar curvature.
\begin{figure}[h]
\begin{picture}(350,110)
\qbezier(50,60)(55,95)(90,100)
\qbezier(50,60)(55,25)(90,20)\qbezier(90,100)(70,60)(90,20)\qbezier[17](90,100)(110,60)(90,20)
\qbezier(90,100)(118,93)(130,75) \qbezier(90,20)(118,27)(130,45)
\qbezier(130,75)(140,70)(180,65) \qbezier(130,45)(140,50)(180,55)
\qbezier(180,65)(220,70)(230,75) \qbezier(180,55)(220,50)(230,45)
\qbezier(240,85)(246,93)(270,100)\qbezier(230,75)(235,82)(240,85)
\qbezier(230,45)(240,27)(270,20) \qbezier(270,100)(305,97)(310,60)
\qbezier(270,20)(305,23)(310,60) \qbezier(40,60)(48,105)(90,110)
\qbezier(40,60)(48,15)(90,10)
\qbezier(90,110)(115,105)(130,85)\qbezier(90,10)(115,15)(130,35)
\qbezier(130,85)(140,80)(205,70) \qbezier(130,35)(140,40)(205,50)
\qbezier(205,70)(237,68)(240,60) \qbezier(205,50)(237,52)(240,60)
\qbezier[35](45,60)(50,100)(90,105)
\qbezier[35](45,60)(50,20)(90,15)
\qbezier[35](90,105)(118,98)(130,80)
\qbezier[33](90,15)(118,22)(130,40)
\qbezier[33](130,80)(140,75)(180,70)
\qbezier[35](130,40)(140,45)(180,50)\qbezier[35](180,70)(220,75)(230,80)
\qbezier[28](180,50)(220,45)(230,40)\qbezier[28](230,80)(235,85)(260,90)
\qbezier[28](230,40)(235,35)(260,30)\qbezier[33](260,90)(285,85)(290,60)
\qbezier[33](260,30)(285,35)(290,60)
\put(300,10){\makebox(0,0)[bl]{$g_1(0)$}}
\put(310,20){\vector(-1,1){10}}
\put(260,5){\makebox(0,0)[bl]{$g_\alpha(0)$}}
\put(260,10){\vector(0,1){20}}
\put(200,20){\makebox(0,0)[bl]{$g_0(0)$}}
\put(210,30){\vector(1,1){25}}
\end{picture}\caption{The smooth family of dumbells}
\end{figure}
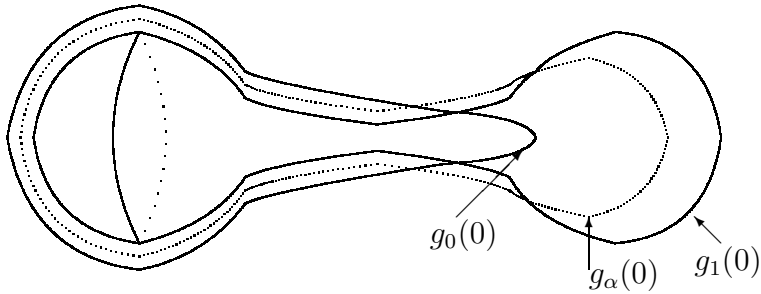

Since the scalar curvature is positive, each solution $g_\alpha(t)$,
$\alpha \in [0,1]$, will exist up to a maximal time
$T_\alpha<+\infty$ and develops a singularity. The main purpose of
this section is to show that there exists $\alpha_0\in[0,1)$ such
that the solution $g_{\alpha_0}(t)$, with the metric
$g_{\alpha_0}(0)$ as initial datum, develops a Type II singularity.
We remark that a Type II singularity might occur in such family of
metrics had been conjectured and the intuition had already described
in \cite{CK} and \cite{T}.

Let us first consider the case that the solutions with degenerate
dumbbells as initial data.

\vskip 0.2cm\noindent {\bf Lemma 4.1}  \emph{ Suppose $g_\alpha (t)$
is a rotationally symmetric solution of the Ricci flow on $S^{n+1}$
with $g_\alpha (0)=\varphi(x,0)^2dx^2+\psi_\alpha ^2(x,0)g_{can}$,
$x\in[-1,1]$. If at initial, the scalar curvature $R^{(\alpha)}>0$,
$\psi_\alpha(x,0)$ has only one bump, it is nonincreasing on the
right polar cap and is strictly concave on the left polar cap, and
$|(\psi_\alpha)_s|(x,0) \leq 1$ on $[-1,1]$, then either the
solution $g_\alpha (t)$ develops a Type II singularity or it shrinks
to a round point.}

\vskip 0.1cm \noindent{\bf Proof.} By the assumption of
$R^{(\alpha)}>0$ and apply the maximum principle to the evolution
equation of the scalar curvature $$\frac{\partial R}{\partial
t}=\Delta R+2|Ric|^2$$ we know that the maximal time $T<+\infty$.

 Now we consider the geometric quantity $s$ defined by
$$s(x,t)=\int_0^x \varphi(x,t)dx.$$ Then the metric can be written as
$$g=ds^2+\psi_\alpha^2(s,t)g_{can}.$$
In the following if we write a relation of the type $f=f(s)$ , it is
to be understood as shorthand for $f=f(s(x,t))$ for evolving
metrics. Since $\psi_\alpha (\pm 1,t)=0$, we know that for any time
$0\leq t<T$, the bump exists. By the standard Sturmian comparison
\cite{An}, we know that $\psi_\alpha(x,t)$ also has a unique bump
for each $t \in [0, T)$. Let $x_\ast (t)$ denote the unique bump. By
Proposition 2.1(3), we can define $$D=\lim_{t\nearrow T}\psi_\alpha
(x_\ast (t),t).$$ We divide it into two cases: \vskip 0.1cm
\noindent{\bf Case 1:} $D>0$.

In this case, by Proposition 2.1(3) and the assumption that
$\psi_\alpha$ is strictly concave on the left polar cap, we know
that no singularity occurs on the left polar cap. Thus the
singularity must occur on the right polar cap. Take the maximal
points $(\tilde{P}_m,t_m)$, i.e., choose the points
$(\tilde{P}_m,t_m)$ such that
$$|Rm(\tilde{P}_m,t_m)|=\sup_{t\leq t_m,Q\in
S^{n+1}}|Rm(Q,t)|\rightarrow+\infty,$$ as $m\rightarrow+\infty$.

Since for any time $t\in[0,T)$ we have $(\psi_\alpha)_s(P,t)=-1$,
where $P$ is the pole of the right polar cap (i.e. the point with
$x=1$), we can choose the nearest point $P_m'$ to $P$ such that
$(\psi_\alpha)_s(P_m',t_m)=-\frac{1}{2}$. If
$d_{t_m}(P_m',P)>d_{t_m}(\tilde{P}_m,P)$, then we set
$P_m=\tilde{P}_m$, otherwise set $P_m=P_m'$. Clearly in the region
between $P_m$ and $P$, we have $|(\psi_\alpha)_s|\geq\frac{1}{2}$.

We first claim that the curvature at $(P_m,t_m)$ is comparable to
the curvature at the maximal point $(\tilde{P}_m,t_m)$. Indeed, if
$P_m=\tilde{P}_m$, then there is nothing to show. If  $P_m\neq
\tilde{P}_m$, then by the estimate in Proposition 2.1(2) and by the
condition that $\psi_\alpha$ is nonincreasing on the right polar cap
and by the choice of the point $P_m$, we know that
$$|Rm(\tilde{P}_m,t_m)|\leq
\frac{C}{\psi_\alpha^2(\tilde{P}_m,t_m)}\leq
\frac{C}{\psi_\alpha^2(P_m,t_m)}.$$ On the other hand
$$
K_1(P_m,t_m)=\frac{1-(\psi_\alpha)_s^2(P_m,t_m)}
{\psi_\alpha^2(P_m,t_m)}=\frac{3}{4\psi_\alpha^2(P_m,t_m)}.
$$ So $$|Rm(P_m,t_m)|\geq K_1(P_m,t_m)=
\frac{3}{4\psi_\alpha^2(P_m,t_m)}\geq
\frac{3}{4C}|Rm(\tilde{P}_m,t_m)|.$$ Obviously since
$(\tilde{P}_m,t_m)$ is the maximal point, we have
$$|Rm(P_m,t_m)|\leq
|Rm(\tilde{P}_m,t_m)|.$$ So the curvatures at $(P_m,t_m)$ and
$(\tilde{P}_m,t_m)$ are comparable, where we used the definition of
$|Rm|$ to be the largest absolute value of the eigenvalues of the
curvature operator $Rm$.

Applying the maximum principle to the evolution equation of the
scalar curvature $R$:$$\frac{\partial R}{\partial t}=\Delta
R+2|Ric|^2$$ and using the pinching estimate in Proposition 2.1(4)
we get
$$\frac{dR_{max}}{dt}\leq CR_{max}^2.$$ Then $$R_{max}(t)\geq
\frac{C}{T-t}$$ for some constant $C$.

We now argue by contradiction to show that the solution develops a
Type II singularity in this case.

Suppose not, then the singularity is of Type I. That is, there
exists some constant $C>0$ such that $$\frac{C^{-1}}{T-t_m}\leq
R^{(\alpha)} (P_m,t_m)\leq \frac{C}{T-t_m}.\eqno(4.1)$$ Define
$$g_{ij}^{(m)}(\cdot,t)=R^{(\alpha)}(P_m,t_m)(g_\alpha)_{ij}(\cdot,t_m+\frac{t}{R^{(\alpha)}(P_m,t_m)}),$$
for $t\in [-t_mR^{(\alpha)}(P_m,t_m),0]$. Then we claim that the
distance from $P_m$ to the pole $P$ measured in the rescaled metric
$g_{ij}^{(m)}(\cdot,0)$ is bounded. Indeed, by the estimate in
Proposition 2.1(2), we know$$|Rm|\leq\frac{C}{\psi_\alpha^2}$$ for
some constant $C$. Then we have
$$\psi_\alpha^2(P_m,t_m)\leq\frac{C}{|Rm(P_m,t_m)|}$$ and using $|(\psi_\alpha)_s|
\geq\frac{1}{2}$ in the region between $P_m$ and $P$, we have
$$d_{t_m}(P_m,P)\leq \frac{\psi_\alpha(P_m,t_m)}{\frac{1}{2}}
\leq\frac{2C}{\sqrt{|Rm(P_m,t_m)|}} \eqno(4.2)$$ where $d_t$ is the
distance measured with the metric $g_\alpha(t)$. Therefore by the
pinching estimate in Proposition 2.1(4) we know that the distance
from $P_m$ to the pole $P$ measured in the rescaled metric
$g_{ij}^{(m)}(\cdot,0)$ is bounded.

 The rescaled
$g_{ij}^{(m)}(t)$ is a solution of the Ricci flow defined for
$t\in[-t_mR^{(\alpha)}(P_m,t_m),0]$ and $0<R^{(m)}(\cdot,t)\leq 1$
and $R^{(m)}(P_m,0)=1$ and has bounded curvature. After taking a
subsequence of $g_{ij}^{(m)}$, we can assume that the marked
manifold $(S^{n+1},g_{ij}^{(m)}(t),P)$ converges to a marked
manifold $(R^{n+1},g_{ij}(t),P), -\infty <t\leq 0$, which is a
solution of the Ricci flow on $R^{n+1}$ with nonnegative curvature
operator ( by the pinching estimate in Proposition 2.1(4) ), has
bounded curvature with $R(P_\ast,0)=1$ at some point $P_\ast$, and
is $\kappa$-noncollapsed for all scales. So the limit is a nonflat
ancient $\kappa$-solution on $R^{n+1}$.

The reduced distance, due to Perelman \cite{P1}, is defined by
$$\arraycolsep=1.5pt\begin{array}{rcl}
&&l^{(\alpha)}(q,\tau)=\frac{1}{2\sqrt{\tau}}\inf\{\int_0^\tau\sqrt{s}
(R^{(\alpha)}(\gamma(s),t_m-s)+|\dot{\gamma}(s)|^2_{(g_\alpha)_{ij}(t_m-s)})ds|\\[4mm]
&&\hskip 3cm\gamma:[0,\tau]\rightarrow S^{n+1}\;\mbox{with}
\;\gamma(0)=P,\gamma(\tau)=q\}.
\end{array}$$
where $\tau=t_m-t, \mbox{for }\;t<t_m $. Then by the Type I
assumption, we have
$$\arraycolsep=1.5pt\begin{array}{rcl}
l^{(\alpha)}(P,\tau)&\leq&\frac{1}{2\sqrt{\tau}}\int_0^\tau\sqrt{s}\frac{C}{T-t_m+s}ds\\[4mm]
&\leq&\frac{C}{2\sqrt{\tau}}\int_0^\tau \frac{1}{\sqrt{s}}ds\\[4mm] &=&
C.
\end{array}\eqno(4.3)$$

We can now use Perelman's backward limit argument in Proposition
11.2 of \cite{P1} to choose a sequence of times $t_k\rightarrow
-\infty$ such that the scaling of $g_{ij}(\cdot,t)$ around $P$ with
the factor $|t_k|^{-1}$ and with the times $t_k$ shifting to the new
time zero converge to a nonflat gradient shrinking soliton in
$C_{loc}^\infty$ topology. Indeed, in the Proposition 11.2 of
\cite{P1}, Perelman takes a limit around some points $q(\tau)$ where
the reduced distance at $q(\tau)$ are uniformly bounded above by
$(n+1)/2$. Instead, in our situation, we want to take a backward
limit around the fixed point $P$. By inspecting the proof of
Proposition 11.2 of \cite{P1} (see also the proof of Theorem 6.2.1
of \cite{CZ} for the details), one only needs to have a uniform
upper bound for the reduced distance at the fixed point $P$. This is
just our estimate (4.3) by the Type I assumption. Then the same
argument as Perelman in section 11.2 in \cite{P1} applies to the
present situation.

By combining with the above Proposition 3.1 and noting that the
gradient shrinking soliton $\bar{g}_{ij}$ is noncompact, we conclude
that the backward limit is $S^n\times R$. But since the limit is
taking around the pole and the metric is rotationally symmetric, it
can not be $S^n\times R$, so we get a contradiction! Hence we have
proved that the singularity is of Type II. \vskip 0.1cm
\noindent{\bf Case 2:} $D=0$.

In this case, if the singularity is of Type II, then there is
nothing to prove. Thus we may assume that the singularity is of Type
I. By the same argument as in the case 1, we can first take a
rescaling limit around the pole $P$ at the maximal time $T$ to get
an ancient $\kappa$-solution and then take a backward limit around
the pole $P$ again to get a nonflat gradient shrinking soliton. If
the shrinking soliton is compact, then by Proposition 3.1 we know
that it is the round $S^{n+1}$. This implies that the original
solution shrinks to a round point as the time tends to the maximal
time $T$. While if the shrinking soliton is noncompact, then by
Proposition 3.1 we know that it is $S^n\times R$; so the same reason
in the proof of the case 1 gives a contradiction! Therefore we have
proved Lemma 4.1.
$$\eqno \#$$

\vskip 0.2cm\noindent {\bf Lemma 4.2} \emph{ The set $A_1$ of
 $\alpha\in[0,1]$ such that the initial metric $g_\alpha(0)$ leads
to a neckpinch singularity of the Ricci flow at some time
$T_\alpha<+\infty$ is open in $[0,1]$.}

\vskip 0.1cm \noindent{\bf Proof.} Obviously it is not empty for
$1\in A_1$.

Suppose $\alpha\in A_1$, then we claim that $\psi_\alpha(0)$ has two
bumps. Otherwise it has only one bump. Since it satisfies the
assumptions of Lemma 4.1 by the above conditions (i), (ii), (iii)
and (iv), then we know that either it shrinks to a round point, or
forms a Type II singularity. Consequently, the initial metric does
not lead to a neckpinch for such $\alpha$. This contradicts with
$\alpha\in A_1$. Similarly, $\psi_\alpha(t)$ has two bumps as long
as the solution exists. Take a small perturbation
$g_\alpha^{(k)}(0)$ of $g_\alpha(0)$ (in $C^3$ topology). Then
$g_\alpha^{(k)}(0)$ still has two bumps. We need to show that
$g_\alpha^{(k)}(0)$ leads to a neckpinch singularity at the maximal
time $T_\alpha^{(k)}<+\infty$.

Since $g_\alpha^{(k)}(0)$ are very close to $g_\alpha(0)$ in $C^3$
topology, the scalar curvatures of the metrics $g_\alpha^{(k)}(0)$
have a uniform positive lower bound. Thus it follows from the
evolution equation of the scalar curvature that the maximal times
$T_\alpha^{(k)}$ are uniformly bounded. After passing to a
subsequence, we can then assume that $T_\alpha^{(k)}\rightarrow
\tilde{T}$ as $k\rightarrow \infty$.

\vskip 0.2cm\noindent {\bf Claim:} $\tilde{T}\geq T$.

Indeed, suppose not, then there exists $\varepsilon >0$, such that
$$\tilde{T}-\varepsilon<T_\alpha^{(k)}<\tilde{T}+\varepsilon<T-\varepsilon<T$$
for all sufficiently large $k$. Consider the time interval
$[0,\tilde{T}-\varepsilon]$. By the assumption that
$g_\alpha^{(k)}(0)$ is sufficiently close to $g(0)$ and $g(t)$ is
smooth on $[0,\tilde{T}-\varepsilon]$, we first show that the
curvature of $g_\alpha^{(k)}(t)$ is uniformly bounded on
$[0,\tilde{T}-\varepsilon]$ for all sufficiently large $k$.

For each $0\leq t \leq \tilde{T}-\varepsilon$, set
$$M(t)=\sup\{|Rm^{(k)}(x,t)| | k\geq 1, x\in S^{n+1}\}$$ and
$$t_0=\sup\{t\geq 0| M(t)<+\infty\},$$ where $Rm^{(k)}$ denotes the
curvature of $g_\alpha^{(k)}(t)$. We want to show that
$t_0=\tilde{T}-\varepsilon$.

By Shi's local derivative estimate in \cite{Sh89}, we know that
$t_0>0$. Suppose $t_0<\tilde{T}-\varepsilon$, then for any small
$\varepsilon'>0$, consider the time interval $[0,t_0-\varepsilon']$.
By the above definition of $M(t)$, we know that the curvature of
$g_\alpha^{(k)}(t)$ is uniformly bounded by $M(t_0-\varepsilon')$ on
$[0,t_0-\varepsilon']$. Take a limit of $g_\alpha^{(k)}(t)$ and by
the uniqueness of the solution to the Ricci flow \cite{Ha82} or
\cite{Ha95F}, we get the limit must be the original solution $g(t)$
on $[0,t_0-\varepsilon']$. So we have the curvature of
$g_\alpha^{(k)}(t)$ is uniformly bounded by some constant $C$ which
does not depend on $\varepsilon'$. Then by Shi's local derivative
estimate in \cite{Sh89} again, we know that the curvature of
$g_\alpha^{(k)}(t)$ is uniformly bounded on
$[0,t_0-\varepsilon'+\frac{1}{C}]$. By choosing $\varepsilon'$ small
enough, we get $t_0-\varepsilon'+\frac{1}{C}>t_0$ and then it is a
contradiction! So we have proved that $t_0=\tilde{T}-\varepsilon$,
that is the curvature of $g_\alpha^{(k)}(t)$ is uniformly bounded on
$[0,\tilde{T}-\varepsilon]$ for all sufficiently large $k$.
Similarly as above, we can take a limit of $g_\alpha^{(k)}(t)$ and
by the uniqueness of the Ricci flow \cite{Ha82} or \cite{Ha95F}, we
get the limit must be $g(t)$ on $[0,\tilde{T}-\varepsilon]$. Using
$g(t)$ is smooth on $[0,\tilde{T}]$, we get the curvature of
$g_\alpha^{(k)}(t)$ is uniformly bounded by some constant $C'$ which
does not depend on $\varepsilon$. Again by Shi's local derivative
estimate in \cite{Sh89}, we know that the curvature of
$g_\alpha^{(k)}(t)$ is uniformly bounded on
$[0,\tilde{T}-\varepsilon+\frac{1}{C'}]$ for all sufficiently large
$k$. Choose $\varepsilon>0$ sufficiently small and using that
$T_\alpha^{(k)}\rightarrow \tilde{T}$ as $k\rightarrow\infty$, we
get $\tilde{T}-\varepsilon+\frac{1}{C'}>T_\alpha^{(k)}$ for all
sufficiently large $k$, which contradicts with the definition of the
maximal time. So $\tilde{T}\geq T$.

Next we show that each $g_\alpha^{(k)}(0)$ leads to a neckpinch
singularity at the maximal time $T_\alpha^{(k)}<+\infty$.

For all sufficiently small $\varepsilon>0$, since $g_\alpha(0)$
leads to a neckpinch, we have
$$\frac{(r_\alpha)_{min}(T-\varepsilon)}
{\psi_\alpha(x_\pm(T-\varepsilon),T-\varepsilon)}\ll 1.\eqno(4.4)$$
By the assumption that $g_\alpha^{(k)}(0)$ is sufficiently close to
$g_\alpha(0)$ and $T_\alpha^{(k)}\rightarrow \tilde{T}\geq T$, we
know that as $k$ large enough, $T_\alpha^{(k)}>T-\varepsilon$ and
$g_\alpha^{(k)}(T-\varepsilon)$ is sufficiently close to
$g_\alpha(T-\varepsilon)$. So
$$\frac{(r_\alpha^{(k)})_{min}(T-\varepsilon)}
{\psi_\alpha^{(k)}(x_\pm(T-\varepsilon),T-\varepsilon)}\ll
1.\eqno(4.5)$$ In views of the work of Angenent and Knopf
\cite{AK}, we know that if we have a rotationally symmetric
$g(0)=ds^2+\psi^2(0)g_{can}$ on $S^{n+1}$ which has two bumps
$x_\pm(0)$ and $\frac{r_{min}(0)}{\psi(x_\pm(0),0)}<C^{-1}$ for
some universal constant $C>0$, (for example we can take $C=100$),
then it leads to a neckpinch singularity. By (4.5), we know that
$g_\alpha^{(k)}(0)$ leads to a neckpinch singularity. Therefore we
proved that $A_1$ is open in $[0,1]$.
$$\eqno \#$$

In the next Proposition, following Perelman's Theorem 12.1 in
\cite{P1}, we will give the singularity structure for the
rotationally symmetric solutions of the Ricci flow. \vskip
0.2cm\noindent {\bf Proposition 4.3} \emph{ Suppose $g_{ij}(t)$,
$t\in[0,T)$, is a rotationally symmetric solution of the Ricci flow
on $S^{n+1}$ with $g(0)=ds^2+\psi^2(0)g_{can}$. If at initial the
scalar curvature $R>0$, then for any given $\varepsilon>0$, there
exists $K=K(\varepsilon, g(0))>\max\{2\varepsilon^{-1},
Q(\frac{3}{4}T)\}>0$, where $Q(\frac{3}{4}T)$ denotes the upper
bound of the curvature  for the times $t\leq\frac{3}{4}T$, such that
for any point $(x_0,t_0)$ with $t_0\geq\frac{3}{4}T$ and
$Q=R(x_0,t_0)\geq K$, the solution in
$\{(y,t)|d_{t_0}^2(y,x_0)<\varepsilon^{-2}Q^{-1},
t_0-\varepsilon^{-2}Q^{-1}\leq t \leq t_0\}$ is, after scaling by
the factor $Q$, $\varepsilon$-close to the corresponding subset of
some orientable ancient $\kappa$-solution, where $\kappa$ is a
positive constant depending only on $T$ and the initial metric
$g(0)$. Consequently, in the region, we have the following gradient
estimates
$$|\nabla(R^{-\frac{1}{2}})|\leq \eta \ \ \mbox{ and } \ \
|\frac{\partial}{\partial t}(R^{-1})|\leq \eta
$$ for some constant $\eta=\eta(\kappa)>0.$} \vskip 0.1cm \noindent{\bf
Proof.} This is just a higher dimensional version of Perelman's
singularity result (Theorem 12.1 of \cite{P1}) for the rotationally
symmetric class. In Theorem 12.1 of \cite{P1}, Perelman obtained
this singularity structure result for any three-dimensional
solution. For the details, one can consult \cite{KL} (from page 83
to  88) or \cite{CZ} (from page 399 to 405). By inspecting
Perelman's argument, when one tries to generalize Perelman's
singularity structure result to higher dimensions, one only needs to
have a higher-dimensional version of the (three-dimensional)
Hamilton-Ivey curvature pinching estimate and shows the canonical
neighborhoods of an ancient $\kappa$-solution consisting
$\varepsilon$-necks or $\varepsilon$-caps. For our case, since the
metric is rotationally symmetric, the estimate due to Angenent-Knopf
in Proposition 2.1(4) gives the desired curvature pinching estimate.
While for a rotationally symmetric ancient $\kappa$-solution, it is
clear that any canonical neighborhood is either an
$\varepsilon$-neck or an $\varepsilon$-cap. So by repeating
Perelman's argument, we obtain the proof of the proposition.
$$\eqno \#$$

We can now prove the main theorem.

\vskip 0.2cm\noindent {\bf Proof of Theorem 1.2.}

Suppose $g_\alpha(t)$ is the family of the solutions to the Ricci
flow satisfies the above conditions (i), (ii), (iii) and (iv). We
want to show that there exists $\alpha_ 0\in [0,1)$ such that for
the solution $g_{\alpha_0}(t)$ of the Ricci flow on $S^{n+1}$ with
initial data $g_{\alpha_ 0}(0)=ds^2+\psi_{\alpha_0}^2g_{can}$,
exists up to a maximal time $T_{\alpha_0}<+\infty$ and develops a
Type II singularity.

Since $g_1(0)=ds^2+\psi_1^2g_{can}$ and by our assumption that
$g_1(0)$ leads to a neckpinch singularity. Then by Lemma 4.2, we
know that $A_1$ is not empty and open in $[0,1]$. While by Lemma
4.1, we know that the solution $g_0(t)$ with the initial data
$g_0(0)$ either develops a Type II singularity or shrinks to a round
point, so $0\overline{\in} A_1$.  Let $(\alpha,1]$ be a connected
component of $A_1$. We want to show that the $\alpha$ is the number
we want.

If $g_\alpha(0)$ develops a Type II singularity, then there is
nothing to show. So in the following we assume it does not develop
a Type II singularity. \vskip 0.1cm \noindent{\bf Claim 1:} \emph{
$\psi_\alpha(0)$ exactly has two bumps.}

Indeed, if $\psi_\alpha(0)$ has only one bump, then $\psi_\alpha$ is
nonincreasing on the right polar cap and by our construction we know
that $\psi_\alpha$ is strictly concave on the left polar cap. So by
Lemma 4.1, we know either the singularity is of Type II or it
shrinks to a round point at the maximal time $T_\alpha<+\infty$. By
our assumption that the singularity is not of Type II. So  it is
shrinking to a round point, and then there exists a time
$\tilde{t}<T_\alpha$ close to $T_\alpha$, such that the curvature is
positive for $t\geq\tilde{t}$. Whenever $\beta\in (\alpha,1]\subset
A_1$ is sufficiently close to $\alpha$, the metric $g_\beta(0)$ is
sufficiently close to the metric $g_\alpha(0)$ ( in the $C^3$
topology). Then by Lemma 4.2 we can choose $\beta\in
(\alpha,1]\subset A_1$ sufficiently close to $\alpha$ so that the
maximal time $T_\beta$ of $g_\beta(t)$ satisfies $T_\beta>\tilde{t}
+ (T_\alpha - \tilde{t})/2$; moreover, by the continuous dependence
of the initial metric, the curvature operator of $g_\beta(t)$ is
also positive at the time $t = \tilde{t}$. Hence by Theorem 1 in
\cite{BW} we know that $g_\beta(0)$ will shrink to a round point at
the maximal time $T_\beta<+\infty$, which contradicts with $\beta\in
(\alpha,1] \subset A_1$. Therefore we have proved the Claim 1.
\vskip 0.1cm \noindent{\bf Claim 2:} \emph{ $\psi_\alpha(t)$ exactly
has two bumps as long as the solution exists.}

Indeed, since $g_\alpha(t)$ is a rotationally symmetric solution of
the Ricci flow on $S^{n+1}$, we know that at the poles of the right
and left polar caps $\psi_\alpha(t)=0$ for any time $0\leq t <
T_\alpha$, so there always exists one bump. By the standard Sturmian
comparison in \cite{An} we know that the number of the bumps is
nonincreasing in time. Suppose at some time $t_0 \in (0,T_\alpha)$
such that the right-most bump disappeared, then $\psi_\alpha(t_0)$
has only one bump. Thus by Lemma 4.1 and our assumption, it shrinks
to a round point. Particularly, there exists a time
$t_0<\tilde{t}<T_\alpha$ such that the curvature is positive for all
times $t\geq \tilde{t}$. By the same argument as above, we can
choose $\beta\in (\alpha,1]\subset A_1$ sufficiently close to
$\alpha$ so that the maximal time $T_\beta$ of $g_\beta(t)$ is
greater than $\tilde{t}$ and the curvature of $g_\beta(t)$ at the
time $t = \tilde{t}$ is also positive. By applying Theorem 1 in
\cite{BW} again we know that $g_\beta(0)$ will shrink to a round
point at the maximal time $T_\beta<+\infty$, which also contradicts
with $\beta\in (\alpha,1]\subset A_1$. Thus we have proved the Claim
2.

So in the following we always assume that $\psi_\alpha(t)$ has two
bumps for all times $t\in [0,T_\alpha)$.

Since at the maximal time $T_\alpha$, the solution $g_\beta(t)$ does
not develop a neckpinch. In views of Angenent-Knopf's result
\cite{AK}, the smaller polar cap must collapse. So, without loss of
generality, we may assume that singularity occurs on the right polar
cap. Similarly as in Lemma 4.1, we first take the maximal points
$(\tilde{P}_m,t_m)$ on $S^{n+1}$ (i.e.,
$|Rm(\tilde{P}_m,t_m)|=\sup_{t\leq t_m,Q\in S^{n+1}}|Rm(Q,t)|$). We
then take the nearest point $P_m'$ to the pole $P$ on the right
polar cap, such that $(\psi_\alpha)_s(P_m',t_m)=-\frac{1}{2}$. If
$d_{t_m}(P_m',P)>d_{t_m}(\tilde{P}_m,P)$, then we set
$P_m=\tilde{P}_m$; otherwise we set $P_m=P_m'$. Clearly in the
region between $P_m$ and $P$, we have
$|(\psi_\alpha)_s|\geq\frac{1}{2}$. \vskip 0.1cm \noindent Define
$$g_{ij}^{(m)}(\cdot,t)=R^{(\alpha)}(P_m,t_m)(g_\alpha)_{ij}(\cdot,t_m+\frac{t}
{R^{(\alpha)}(P_m,t_m)})$$ for $t\in
[-t_mR^{(\alpha)}(P_m,t_m),0]$.

\vskip 0.1cm \noindent{\bf Claim 3:}  \emph{ A subsequence of
$g_{ij}^{(m)}(\cdot,t)$ around the point $P$ will converge to a
nonflat complete ancient $\kappa$-solution on a smooth manifold $M$,
where $\kappa$ is some positive constant depending only on the
initial metric $g_\alpha(0)$. }

Indeed, if the maximal point $(\tilde{P}_m,t_m)$ is on the right
polar cap, then by the estimate in Proposition 2.1(2) and by the
condition that $\psi_\alpha$ is nonincreasing on the right polar cap
and by the choice of the point $P_m$, we know that
$$|Rm(\tilde{P}_m,t_m)|\leq
\frac{C}{\psi_\alpha^2(\tilde{P}_m,t_m)}\leq
\frac{C}{\psi_\alpha^2(P_m,t_m)}.$$ On the other hand
$$
K_1(P_m,t_m)=\frac{1-(\psi_\alpha)_s^2(P_m,t_m)}
{\psi_\alpha^2(P_m,t_m)}=\frac{3}{4\psi_\alpha^2(P_m,t_m)}.
$$ So $$|Rm(P_m,t_m)|\geq K_1(P_m,t_m)=
\frac{3}{4\psi_\alpha^2(P_m,t_m)}\geq
\frac{3}{4C}|Rm(\tilde{P}_m,t_m)|.$$ Obviously since
$(\tilde{P}_m,t_m)$ is the maximal point, we have
$$|Rm(P_m,t_m)|\leq
|Rm(\tilde{P}_m,t_m)|.$$ So the curvatures at $(P_m,t_m)$ and
$(\tilde{P}_m,t_m)$ are comparable, where we used the definition of
$|Rm|$ to be the largest absolute value of the eigenvalues of the
curvature operator $Rm$. So by repeating (part of) the argument as
in case 1 in Lemma 4.1, we know that a subsequence of
$g_{ij}^{(m)}(\cdot,t)$ around the point $P$ will converge to a
nonflat complete ancient $\kappa$-solution on a smooth manifold $M$
for some positive constant $\kappa$ depending only on the initial
metric $g_\alpha(0)$.

We remain to consider the case that the maximal point
$(\tilde{P}_m,t_m)$ does not lie on the right polar cap, then it
must lie in the region between the two bumps.

In this case, we first prove the following assertion:

\emph{ For any $A<+\infty$, there exists a positive constant $C(A)$
such that the curvatures of $g_{ij}^{(m)}(\cdot,t)$ at the new time
$t=0$ satisfy the estimate$$|Rm^{(m)}(y,0)|\leq C(A)$$ whenever
$d_{g^{(m)}(\cdot,0)}(y,P_m)\leq A$ and $m\geq 1$, where $Rm^{(m)}$
denotes the curvature of the metric $g_{ij}^{(m)}$.}

This assertion in the three-dimensional Ricci flow has been verified
by Perelman in his proof of the Theorem 12.1 in \cite{P1} (the first
detailed exposition of this part of Perelman's argument appeared in
the first version of Kleiner-Lott \cite{KL}), where the only
three-dimension features he used are the Hamliton-Ivey curvature
pinching estimate and the canonical neighborhood condition of an
ancient $\kappa$-solution consisting the $\varepsilon$-necks and
$\varepsilon$-caps. In our case, by noting that the metric is
rotationally symmetric, the canonical neighborhood condition can be
easily obtained as pointed out before, and the pinching estimate has
already given in Proposition 2.1(4). So by some slight
modifications, Perelman's argument also works for our case. In the
following we only give the details for the modified parts. For the
complete details, one can compare with \cite{KL} (from page 85 to
87) or \cite{CZ} (from page 400 to 402).

For each $\rho \geq 0$, set
$$M(\rho)=\sup\{R^{(m)}(x,0)\ |\ m\geq 1, x\in S^{n+1} \ \mbox{ with } \
d_{0}(x,P_m)\leq \rho\}
$$
and
$$\rho_0=\sup\{\rho \geq 0\  |\ \ M(\rho)<+\infty \}.
$$
By the pinching estimate in Proposition 2.1(4), it suffices to show
$\rho_0=+\infty.$

We need to adapt Perelman's argument to show that $\rho_0>0$.

For arbitrary fixed small $\varepsilon >0$, by Proposition 4.3, we
know that there exists $K=K(\varepsilon,
g_\alpha(0))>\max\{2\varepsilon^{-1}, Q(\frac{3}{4}T_\alpha)\}>0$,
where $Q(\frac{3}{4}T_\alpha)$ denotes the upper bound of the
curvature for the times $t\leq\frac{3}{4}T_\alpha$, such that for
any point $(x_0,t_0)$ with $t_0\geq\frac{3}{4}T_\alpha$ and
$Q=R^{(\alpha)}(x_0,t_0)\geq K$, the solution in
$\{(y,t)|d_{t_0}^2(y,x_0)<\varepsilon^{-2}Q^{-1},
t_0-\varepsilon^{-2}Q^{-1}\leq t \leq t_0\}$ is, after scaling by
the factor $Q$, $\varepsilon$-close to the corresponding subset of
some orientable ancient $\kappa$-solution for some positive constant
$\kappa$ depending only on the initial metric $g_\alpha(0)$.
Consequently we have the gradient estimate in the
region$$|\nabla(R^{-\frac{1}{2}})|\leq \eta \ \ \mbox{ and } \ \
|\frac{\partial}{\partial t}(R^{-1})|\leq \eta \eqno(4.6)$$ for some
constant $\eta=\eta(\kappa)>0.$

If $R^{(\alpha)}(P_m,t_m)\geq K$, then by the above gradient
estimate (4.6), we know that there exists some constant
$c=c(\eta)>0$ such that $$R^{(\alpha)}(x,t_m)\leq
2R^{(\alpha)}(P_m,t_m)$$ for any point $x\in
B_{t_m}(P_m,c(R^{(\alpha)}(P_m,t_m))^{-\frac{1}{2}})$. Hence in
this case we have $\rho_0\geq c>0$.

If $R^{(\alpha)}(P_m,t_m)< K$, then we prove that $\rho_0\geq
\tilde{c}$ for some constant $\tilde{c}=\tilde{c}(c,K,\bar{c})$,
where $\bar{c}$ is the positive lower bound of the scalar
curvature $R^{(\alpha)}$ on $S^{n+1}$ at initial time. In fact,
consider the points $x\in
B_{t_m}(P_m,\frac{c}{2}(R^{(\alpha)}(P_m,t_m))^{-\frac{1}{2}})$,
if $R^{(\alpha)}(x,t_m)< K$ for all points $x$, then $\rho_0\geq
\frac{c}{2}>0$ (since $R^{(\alpha)}(P_m,t_m)\geq \bar{c})$; if
$R^{(\alpha)}(x,t_m)\geq K$ for some point $x$, consider the
nearest point $y_0 \in
B_{t_m}(P_m,\frac{c}{2}(R^{(\alpha)}(P_m,t_m))^{-\frac{1}{2}})$ to
$P_m$ such that $R^{(\alpha)}(y_0,t_m)= K$, then by Proposition
4.3 and the gradient estimate (4.6), we know that
$$R^{(\alpha)}(y,t_m)\leq 2R^{(\alpha)}(y_0,t_m)=2K \eqno(4.7)$$
for any point $y\in
B_{t_m}(y_0,c(R^{(\alpha)}(y_0,t_m))^{-\frac{1}{2}})$. Since the
scalar curvature has a positive lower bound $\bar{c}$ by our
assumption, we know that there exists
$c'=cK^{-\frac{1}{2}}\bar{c}^{\frac{1}{2}}>0$ such that
$$B_{t_m}(y_0,c(K)^{-\frac{1}{2}})\supset
B_{t_m}(P_m,c'(R^{(\alpha)}(P_m,t_m))^{-\frac{1}{2}}).\eqno(4.8)$$
By (4.7) and (4.8) we know that for any point $z\in
B_{t_m}(P_m,c'(R^{(\alpha)}(P_m,t_m))^{-\frac{1}{2}})$ we have
$$R^{(\alpha)}\leq 2K.$$ Then $\rho_0\geq c'>0$.
Set $\tilde{c}=\min\{\frac{c}{2},c'\}$, then in this case we have
$\rho_0\geq \tilde{c}>0$.

Hence we have proved $\rho_0>0$.

In the rest, we can apply the same argument of Perelman \cite{P1}
(see also \cite{KL} and \cite{CZ} for details) to obtain that
$\rho_0=+\infty$. That is, the curvatures of $g_{ij}^{(m)}(\cdot,t)$
at the new times $t=0$ stay uniformly bounded at bounded distances
from $P_m$ for all $m$. Furthermore, by the estimate in Proposition
2.1(2) and using $|(\psi_\alpha)_s|\geq\frac{1}{2}$ in the region
between $P_m$ and $P$, we know that the distance from $P_m$ and $P$
measured in the rescaled metric $g_{ij}^{(m)}(\cdot,0)$ is bounded.
So we obtained that the curvature of $g_{ij}^{(m)}(\cdot,t)$ at the
new times $t=0$ stay uniformly bounded at bounded distances from $P$
for all $m$. This completes the proof of the assertion.

By the gradient estimate in Proposition 4.3 and Shi's local
derivative estimate in \cite{Sh89} and Hamilton's compactness
theorem in \cite{Ha95}, we can take a $C_{loc}^\infty$ subsequent
limit to obtain $(M,g_\infty(\cdot,t),P)$ which is complete,
$\kappa$-noncollapsed on all scales and is defined on a space-time
open subset of $M\times (-\infty,0]$ containing the time slice
$M\times \{0\}$. Clearly it follows from the pinching estimate in
Proposition 2.1(4) that the limit $(M,g_\infty(\cdot,t),P)$ has
nonnegative curvature operator. Then exactly as Perelman's argument
in Theorem 12.1 of \cite{P1} (see also \cite{KL} and \cite{CZ} for
details), we can get that the curvature of the limit
$g_\infty(\cdot,t)$ at $t=0$ has bounded curvature and also that the
limit $g_\infty(\cdot,t)$ can be defined on $(-\infty,0]$. So we
have proved that $g_\infty(\cdot,t)$ is an ancient $\kappa$-solution
on $M$ and Claim 3 holds.

Since by our assumption that the singularity is not of Type II. Then
there exists some constant $\tilde{C}>0$ such that
$$0\leq R(P_m,t)\leq
\frac{\tilde{C}}{T_\alpha -t}.$$ Then by Claim 3 we know that a
subsequence of $g_{ij}^{(m)}(\cdot,t)$ around $P$ converges to a
nonflat ancient $\kappa$-solution $g_{ij}$ on $M$.  Then by the same
proof as in the case 1 in Lemma 4.1, we obtain that there exists a
sequence of times $t_k\rightarrow -\infty$ such that the scaling of
$g_{ij}(\cdot,t)$ around $P$ with the factor $|t_k|^{-1}$ and with
the times $t_k$ shifting to the new time zero converge to a nonflat
gradient shrinking soliton in $C_{loc}^\infty$ topology. If the
nonflat gradient shrinking soliton is noncompact, Proposition 3.1
gives us that it is $R\times S^n$. But since the limit is taking
around the pole $P$ and the metric is rotationally symmetric, it can
not be $R\times S^n$. So this contradiction implies that the nonflat
gradient shrinking soliton is compact. By Proposition 3.1 again, we
know that it is the round $S^{n+1}$. Consequently the curvature of
the original solution becomes positive as the time $t$ close to the
maximal time $T_\alpha$. Then repeat the same proof as in Claim 1,
we can choose $\beta\in (\alpha,1]\subset A_1$ sufficiently close to
$\alpha$ such that $g_\beta(0)$ will also shrink to a round point at
the maximal time $T_\beta<+\infty$, which contradicts with $\beta\in
(\alpha,1] \subset A_1$. So the singularity must be of Type II.
Therefore we have proved our theorem 1.2.$$\eqno \#$$

\vskip 0.2cm\noindent {\bf Remark 1.} During the proof of the main
theorem, we actually proved the existence of Type II singularities
on noncompact manifolds. More precisely, we proved that for each
$n\geq2$, there exists complete and rotationally symmetric metrics
on $R^{n+1}$ with bounded curvatures such that the Ricci flow
starting at the metrics develop Type II singularities at some times
$T<+\infty$. In particular, we can take the initial metrics on
$R^{n+1}$ to be the complete and rotationally symmetric, with
nonnegative sectional curvature and positive scalar curvature, and
 asymptotic to the round cylinder of scalar curvature 1 at infinity. \vskip
0.2cm\noindent {\bf Remark 2.} In the unpublished preprint
\cite{BR}, Robert Bryant proved the existence of the nontrivial
steady Ricci solitons on $R^{n}$ by solving certain nonlinear ODE
system. These steady Ricci solitons are complete, rotationally
symmetric with positive curvatures. By combining with the work of
Hamilton \cite{Ha93}, this paper gives another proof for the
existence of the nontrivial steady Ricci solitons on $R^n$ for all
dimensions $n\geq 3$, which are also complete, rotationally
symmetric and have positive curvatures.

\end{document}